\documentclass[10pt,leqno]{article}

\usepackage{amsmath,amssymb,amsthm,mathrsfs,dsfont}

\usepackage[margin=3cm]{geometry} 

\usepackage{titlesec,hyperref}

\usepackage{color}

\usepackage{fancyhdr}
\titleformat{\subsection}{\it}{\thesubsection.\enspace}{1pt}{}

\newtheorem{theo}{Theorem}[section]
\newtheorem{lemm}[theo]{Lemma}

\newtheorem{prop}[theo]{Proposition}
\newtheorem{rema}[theo]{Remark}
\numberwithin{equation}{section}

\allowdisplaybreaks 

\begin{document}
\title{Ill-posedness for the Cauchy problem of the modified Camassa-Holm equation in $B_{\infty,1}^0$
\hspace{-4mm}
}

\author{Zhen $\mbox{He}^1$ \footnote{Email: hezh56@mail2.sysu.edu.cn},\quad
	Zhaoyang $\mbox{Yin}^{1,2}$\footnote{E-mail: mcsyzy@mail.sysu.edu.cn}\\
	$^1\mbox{Department}$ of Mathematics,
	Sun Yat-sen University, Guangzhou 510275, China\\
	$^2\mbox{School}$ of Science,\\ Shenzhen Campus of Sun Yat-sen University, Shenzhen 518107, China}

\date{}
\maketitle
\hrule

\begin{abstract}
In this paper, we prove the norm inflation and get the ill-posedness for the modified Camassa-Holm equation in $B_{\infty,1}^0$. Therefore we completed all well-posedness and ill-posedness problem for the modified Camassa-Holm equation in all critical spaces $B_{p,1}^\frac{1}{p}$ with $p\in[1,\infty]$.   \\
\vspace*{5pt}

\vspace*{5pt}
\end{abstract}

\vspace*{10pt}

\tableofcontents

\section{Introduction}
 As far as we know, the CH equation has many properties, such as: integrability  \cite{Constantin2001,Constantin2006,Constantin1999},  Hamiltonian structure, infinitely many conservation laws\cite{Fuchssteiner1981,Constantin2001}.  The local well-posedness  of  CH equation in Sobolev spaces $H^s, s>\frac 3 2$ and in Besov space $B^{s}_{p,r}$ with $ s>\max\{\frac{3}{2},1+\frac 1 p \}$ or $s=1+\frac 1 p,p\in [1,2],r=1$ has benn studied in \cite{Constantin1998,Constantin1998w-l,Danchin2001inte,Danchin2003wp,Liu2019ill,Himonas2012CN,Li2016nwpC}. In \cite{Constantin2000gw,Xin2000ws}, the authors showed  the existence and uniqueness of global weak solutions for the CH equation. Later, Bressan and Constantin turned the equation to a semilinear system to studied the existence of the global conservative solutions \cite{book}  and global disspative solutions \cite{Bressan2007gd}  in $H^1(\mathbb{R}).$  Recently, Ye, Yin and Guo \cite{2021} gave the local well-posedness in $B^{1+\frac 1 p}_{p,1}$ with $p\in [1,+\infty].$ For the ill-posed problems,  the author \cite{Byers2006illhs} proved that CH equation is ill-posed in $H^{s}(\mathbb{R}), s<\frac 3 2.$ Danchin  proved that the CH equation is ill-posed in $B^{\frac 3 2}_{p,\infty} $  \cite{Danchin2001inte,Danchin2003wp}  and Guo et al. discussed  the equation for the shallow water wave type is ill-posed in $H^{\frac 3 2}$ and  $B^{1+\frac 1 p}_{p,r}$ with $ p\in [1,+\infty], r\in (1,+\infty]$   \cite{Liu2019ill}. In \cite{Li2021}, Li, Yu et.al  studied the ill-posedness of the Camassa-Holm equation in $B^{s}_{p,\infty}$ for $s>2+\max\{\frac 3 2, 1+\frac 1  p\}, 1\leq p\leq +\infty.$
 
     ~~ If we study the Camassa-Holm through the geometric approach in\cite{Chern} \cite{Reyes}, we will obtain the following modified Camassa-Holm (MOCH) equation,
   \begin{equation}\label{eq0}
   	\left\{\begin{array}{l}
   		\gamma_t=\lambda(v_x-\gamma-\frac{1}{\lambda}v\gamma)_x
   		,  \quad t>0,\ x\in\mathbb{R},  \\
   		v_{xx}=\gamma_x+\frac {\gamma^2} {2\lambda} ,  \quad t\geq0,\ x\in\mathbb{R},  \\
   		\gamma(0,x)=\gamma_0(x),  \quad x\in\mathbb{R},
   	\end{array}\right.
   \end{equation}
   which was called by Gorka and Reyes\cite{GR2010}. Let $G={\partial_x}^2-1,m=Gv$.
   
   The equation \eqref{eq0} can be rewritten as
   \begin{equation}\label{eq1}
   	\left\{\begin{array}{l}
   		\gamma_t+G^{-1}m\gamma_x=\frac {\gamma^2}{2} +\lambda G^{-1}m-
   		\gamma G^{-1}m_x,  \quad t>0,\ x\in\mathbb{R},  \\
   		m=\gamma_x+\frac {\gamma^2} {2\lambda} ,  \quad t\geq0,\ x\in\mathbb{R},  \\
   		\gamma(0,x)=\gamma_0(x),  \quad x\in\mathbb{R}.
   	\end{array}\right.
   \end{equation}

   Conservation laws and the existence and uniqueness of weak solutions to the
   modified Camassa-Holm equation were presented in \cite{GR2010}. We observe that if we solve
   \eqref{eq1}, then $m$ will formally satisfy the following physical form of the Camassa-Holm
   
   \begin{align}
   	m_t=-2vm_x-mv_x+\lambda v_x
   \end{align}

   If $\lambda=0$, it is known as the well-known Camassa-Holm (CH) equation. 
   Luo, Qiao and Yin studied the locally well-posedness in $B^s_{p,r},s\textgreater \max\{\frac{1}{2},\frac{1}{p}\}$ or $s=\frac{1}{p},1\leq p\leq 2,r=1$,blow up condition, global existence for periodic MOCH and global conservative solution\cite{Luo1,Luo2,Luo3}. And the authors constructed local well-posednesss for the Cauchy problem of a modified Camassa-Holm (MOCH) equation in nonhomogeneous Besov spaces $B^{\frac 1 p}_{p,1}$ with $1\leq p<+\infty$, which completed the local well-posedness problem in Besov space in \cite{He}.
   
   This paper is aimed to prove the locally ill-posedness in critical Besov space $B^{0}_{\infty,\infty,1}$. Motivated by \cite{Yei}, we refine the space $B^{0}_{\infty,1}$ with the norm $\|u\|_{B_{\infty,\infty,1}^{0}}=\sup_j (j+2)^2\|\Delta_{j}u\|$, which can help us to get the information of the low frequencies of the function while the good properties found in \cite{Yei} still holds ture.
   The content of this paper is the following. In Section 2, we recall some basic definitions and the
   related results about Besov spaces. Section 3 presents some detailed proof of some calculations we will use later. In Section 4, we study about local well-posedness of the modified Camassa-Holm in $C\Big([0,T];B^{0}_{\infty,1}(\mathbb{R})\cap B^{0}_{\infty,\infty,1}(\mathbb{R})\Big)\cap C^1\Big([0,T];B^{-1}_{\infty,1}(\mathbb{R})\cap B^{-1}_{\infty,\infty,1}(\mathbb{R})\Big)$. Section 5 devotes
   to showing the norm inflation to prove the ill-posedness.
   
   We state our two main theorems as follows:
   \begin{theo}\label{them2}
   	Let $\gamma^0 \in B^{0}_{\infty,1}(\mathbb{R})\cap B^{0}_{\infty,\infty,1}(\mathbb{R})$.
   	Then there exists a time $T>0$ such that \eqref{eq1} has a unique solution $\gamma \in E^p_T\triangleq C\Big([0,T];B^{0}_{\infty,1}(\mathbb{R})\cap B^{0}_{\infty,\infty,1}(\mathbb{R})\Big)\cap C^1\Big([0,T];B^{-1}_{\infty,1}(\mathbb{R})\cap B^{-1}_{\infty,\infty,1}(\mathbb{R})\Big).$  Moreover,  the solution depends continuously on the initial data.
   \end{theo}
\begin{theo}\label{the1}
	For any 10 \textless N $\in \mathbb{N}^+$ large enough, there exists a $u_0\in \mathcal{C}^\infty(\mathbb{R})$ such that the  following hold:
	\\
	1. $\|\gamma_0\|_{B_{\infty,1}^0}\leq CN^{-\frac{1}{10}}$
	\\
	2. There is a unique solution $\gamma \in \mathcal{C}_T(\mathcal{C}^\infty(R))$ to the Cauchy problem \eqref{eq1} with a time T$\leq N^{-\frac{1}{2}}$
	\\
	3. There exists a time $t_0\in [0,T]$ such that $\|\gamma(t_0)\|_{B_{\infty,1}^0}\geq \ln N$
\end{theo}

\section{Preliminaries}
  ~~~~In this section, we will present some propositions about the Littlewood-Paley decomposition and Besov spaces.
  \begin{prop}\cite{Chemin2011,He2017}\label{Bernstein}
  	Let $\mathscr{B}$ be a ball and $\mathscr{C}$ be an annulus. A constant $C\textgreater0$ exists such that for all $k\in\mathbb{N}$, $1\leq p\leq q\leq \infty$, and any function $f\in L^p(\mathbb{R})$, we have \\
  	$$Supp(\mathscr {F}f) \subset \lambda\mathscr{B} \Longrightarrow \|D^k f\|_{L^p}=\sup\limits_{|\alpha|=k}\|\partial^{\alpha}f\|_{L^q} \leq C^{k+1}{\lambda}^{k+d(\frac{1}{p}-\frac{1}{q})}\|f\|_{L^p}$$
  		$$Supp(\mathscr {F}f) \subset \lambda\mathscr{C} \Longrightarrow  C^{-k-1}{\lambda}^{k}\|f\|_{L^p} \leq \|D^k f\|_{L^p} \leq C^{k+1}{\lambda}^{k}\|f\|_{L^p}$$
  \end{prop}
  We give two useful interpolation inequalities.
  \begin{prop}\label{prop}\cite{Chemin2011,He2017}
  		Let $\mathcal{C}$ be the annulus $\{\xi\in\mathbb{R}^d:\frac 3 4\leq|\xi|\leq\frac 8 3\}$. There exist radial function $\varphi$, valued in the interval $[0,1]$, belonging respectively to $\mathcal{D}(\mathcal{C})$, and such that
  	$$ \forall\xi\in\mathbb{R}^d\backslash\{0\},\ \sum_{j\in\mathbb{Z}}\varphi(2^{-j}\xi)=1, $$
  	$$ |j-j'|\geq 2\Rightarrow\mathrm{Supp}\ \varphi(2^{-j}\cdot)\cap \mathrm{Supp}\ \varphi(2^{-j'}\cdot)=\emptyset. $$
  	Further, we have
  	$$ \forall\xi\in\mathbb{R}^d\backslash\{0\},\ \frac 1 2\leq\sum_{j\in\mathbb{Z}}\varphi^2(2^{-j}\xi)\leq 1. $$

  	Let $u$ be a tempered distribution in $\mathcal{S}'_h(\mathbb{R}^d)$. For all $j\in\mathbb{Z}$, define
  	$$
  	\dot{\Delta}_j u=\mathcal{F}^{-1}(\varphi(2^{-j}\cdot)\mathcal{F}u).
  	$$
  	Then the Littlewood-Paley decomposition is given as follows:
  	$$ u=\sum_{j\in\mathbb{Z}}\dot{\Delta}_j u \quad \text{in}\ \mathcal{S}'(\mathbb{R}^d). $$
  	Now, we introduce the definition of homogenous Besov spaces as follows.
  \end{prop}
  	\begin{prop}\cite{Chemin2011}
  	Let $s\in\mathbb{R},\ 1\leq p,r\leq\infty.$
  	\begin{equation*}\left\{
  		\begin{array}{l}
  			B^s_{p,r}\times B^{-s}_{p',r'}\longrightarrow\mathbb{R},  \\
  			(u,\phi)\longmapsto \sum\limits_{|j-j'|\leq 1}\langle \Delta_j u,\Delta_{j'}\phi\rangle,
  		\end{array}\right.
  	\end{equation*}
  	defines a continuous bilinear functional on $B^s_{p,r}\times B^{-s}_{p',r'}$. Denoted by $Q^{-s}_{p',r'}$ the set of functions $\phi$ in $\mathcal{S}'$ such that
  	$\|\phi\|_{B^{-s}_{p',r'}}\leq 1$. If $u$ is in $\mathcal{S}'$, then we have
  	$$\|u\|_{B^s_{p,r}}\leq C\sup_{\phi\in Q^{-s}_{p',r'}}\langle u,\phi\rangle.$$
  \end{prop}
  We then have the following product laws:
  \begin{lemm}\label{product}\cite{Chemin2011,He2017}
  	(1) For any $s>0$ and any $(p,r)$ in $[1,\infty]^2$, the space $L^{\infty} \cap B^s_{p,r}$ is an algebra, and a constant $C=C(s,d)$ exists such that
  	$$ \|uv\|_{B^s_{p,r}}\leq C(\|u\|_{L^{\infty}}\|v\|_{B^s_{p,r}}+\|u\|_{B^s_{p,r}}\|v\|_{L^{\infty}}). $$
  	(2) If $1\leq p,r\leq \infty,\ s_1\leq s_2,\ s_2>\frac{d}{p} (s_2 \geq \frac{d}{p}\ \text{if}\ r=1)$ and $s_1+s_2>\max(0, \frac{2d}{p}-d)$, there exists $C=C(s_1,s_2,p,r,d)$ such that
  	$$ \|uv\|_{B^{s_1}_{p,r}}\leq C\|u\|_{B^{s_1}_{p,r}}\|v\|_{B^{s_2}_{p,r}}. $$
  	(3) If $1\leq p\leq 2$,  there exists $C=C(p,d)$ such that
  	$$ \|uv\|_{B^{\frac d p-d}_{p,\infty}}\leq C \|u\|_{B^{\frac d p-d}_{p,\infty}}\|v\|_{B^{\frac d p}_{p,1}}. $$
  \end{lemm}
  The Gronwall lemma as follows.
  \begin{lemm}\label{osgood}\cite{Chemin2011}
  	Let $f(t),~ g(t)\in C^{1}([0,T]), f(t), g(t)>0.$ Let $h(t)$ is a continuous function on $[0,T].$ Assume that, for any $t\in [0,T]$ such that
  	$$\frac 1 2 \frac{d}{dt}f^{2}(t)\leq h(t)f^{2}(t)+g(t)f(t).$$
  	Then for any time $t\in [0,T],$ we have
  	$$f(t)\leq f(0)exp\int_0^th(\tau)d\tau+\int_0^t g(\tau)\ exp(\int_\tau ^t h(\tau)dt')d\tau.$$
  \end{lemm}
  Now we state some useful results in the transport equation theory, which are important to the proofs of our main theorem later.
  \begin{equation}\label{transport}
  	\left\{\begin{array}{l}
  		f_t+v\cdot\nabla f=g,\ x\in\mathbb{R}^d,\ t>0, \\
  		f(0,x)=f_0(x).
  	\end{array}\right.
  \end{equation}
  \begin{lemm}\label{existence}\cite{Chemin2011}
  	Let $1\leq p\leq p_1\leq\infty,\ 1\leq r\leq\infty,\ s> -d\min(\frac 1 {p_1}, \frac 1 {p'})$. Let $f_0\in B^s_{p,r}$, $g\in L^1([0,T];B^s_{p,r})$, and let $v$ be a time-dependent vector field such that $v\in L^\rho([0,T];B^{-M}_{\infty,\infty})$ for some $\rho>1$ and $M>0$, and
  	$$
  	\begin{array}{ll}
  		\nabla v\in L^1([0,T];B^{\frac d {p_1}}_{p_1,\infty}), &\ \text{if}\ s<1+\frac d {p_1}, \\
  		\nabla v\in L^1([0,T];B^{s-1}_{p,r}), &\ \text{if}\ s>1+\frac d {p_1}\ or\ (s=1+\frac d {p_1}\ and\ r=1).
  	\end{array}
  	$$
  	Then the equation \eqref{transport} has a unique solution $f$ in \\
  	-the space $C([0,T];B^s_{p,r})$, if $r<\infty$; \\
  	-the space $\Big(\bigcap_{s'<s}C([0,T];B^{s'}_{p,\infty})\Big)\bigcap C_w([0,T];B^s_{p,\infty})$, if $r=\infty$.
  \end{lemm}
  \begin{lemm}\label{priori estimate}\cite{Chemin2011,Li2016nwpC}
  	Let $s\in\mathbb{R},\ 1\leq p,r\leq\infty$.
  	There exists a constant $C$ such that for all solutions $f\in L^{\infty}([0,T];B^s_{p,r})$ of \eqref{transport} in one dimension with initial data $f_0\in B^s_{p,r}$, and $g\in L^1([0,T];B^s_{p,r})$, we have for a.e. $t\in[0,T]$,
  	$$ \|f(t)\|_{B^s_{p,r}}\leq \|f_0\|_{B^s_{p,r}}+\int_0^t \|g(t')\|_{B^s_{p,r}}dt'+\int_0^t V^{'} (t^{'})\|f(t)\|_{B^s_{p,r}}dt{'} $$
  	or
  	$$ \|f(t)\|_{B^s_{p,r}}\leq e^{CV(t)}\Big(\|f_0\|_{B^s_{p,r}}+\int_0^t e^{-CV(t')}\|g(t')\|_{B^s_{p,r}}dt'\Big) $$
  	with
  	\begin{equation*}
  		V'(t)=\left\{\begin{array}{ll}
  			\|\nabla v\|_{B^{s+1}_{p,r}},\ &\text{if}\ s>\max(-\frac 1 2,\frac 1 {p}-1), \\
  			\|\nabla v\|_{B^{s}_{p,r}},\ &\text{if}\ s>\frac 1 {p}\ \text{or}\ (s=\frac 1 {p},\ p<\infty, \ r=1),
  		\end{array}\right.
  	\end{equation*}
  	and when $s=\frac 1 p-1,\ 1\leq p\leq 2,\ r=\infty,\ \text{and}\ V'(t)=\|\nabla v\|_{B^{\frac 1 p}_{p,1}}$.\\
  	If $f=v,$ for all $s>0, V^{'}(t)=\|\nabla v(t)\|_{L^{\infty}}.$
  \end{lemm}

  \begin{lemm}\label{cont1}\cite{Chemin2011,Li2016nwpC}
  	Let $y_0\in B^{\frac{1}{p}}_{p,1}$ with $1\leq p<\infty,$  and $f\in L^1([0,T];B^{\frac{1}{p}}_{p,1}).$ Define $\bar{\mathbb{N}}=\mathbb{N} \cup {\infty},$ for $n\in \bar{\mathbb{N}},$ denote by $y_n \in C([0,T];B^{\frac{1}{p}}_{p,1})$ the solution of
  	\begin{equation}
  		\left\{\begin{aligned}
  			&\partial_ty_n+A_n(u)\partial_xy_n=f,\quad x\in \mathbb{R},\\
  			&y_n(t,x)|_{t=0}=y_0(x). \\
  		\end{aligned} \right. \label{e1}
  	\end{equation}
  	Assume for some $\alpha(t)\in L^1(0,T),\  \sup\limits_{n\in \bar{\mathbb{N}}} \|A_n(u)\|_{B^{1+\frac{1}{p}}_{p,1}}\leq \alpha (t).$ If $A_n(u)$ converges in $A_{\infty}(u)$ in $L^1([0,T];B^{\frac{1}{p}}_{p,1}),$ then the sequence $(y_n)_{n\in \mathbb{N}}$ converges in $C([0,T];B^{\frac{1}{p}}_{p,1}).$
  \end{lemm}
  Let us consider the following initial value problem
  \begin{equation}\label{eq101}
  	\left\{\begin{array}{l}
  		y_t=u(t,y),t\in [0,T)  \\
  		
  		y(0,x)=x,x\in\mathbb{R}.
  	\end{array}\right.
  \end{equation}
  \begin{lemm}\cite{Yinper}\cite{Conper}
  	Let u $\in C( [ 0,T);H^s)\cap  C^1([0,T);H^{s-1}),s\geq 2.$Then\eqref{eq101} has a unique solution $y \in C^1([0,T)\times \mathbb{R};\mathbb{R})$.Moreover,the map $y(t,\cdot)$ is an increasing diffeomorphism of $\mathbb{R}$ with  
  	$$y_x(t,x)=exp(\int_{0}^{t} u_x(s,q(s,x))ds)\textgreater 0,\forall(t,x)\in [0,T)\times \mathbb{R}$$ 
  \end{lemm}
   \begin{lemm}\label{2continuity}\cite{Li2016nwpC}
  	Suppose that ~$1\leq p\leq\infty,\ 1\leq r<\infty,\ s>\frac d p$\ (or \ $s=\frac d p,\ p<\infty,\ r=1)$. Let ~$\bar{\mathbb{N}}=\mathbb{N}\cup\{\infty\}$. Asuming that ~$(v^n)_{n\in\bar{\mathbb{N}}}\in C([0,T];B^{s+1}_{p,r})$, and ~$(f^n)_{n\in\bar{\mathbb{N}}} \in C([0,T];B^s_{p,r})$ solves the equation:
  	\begin{equation}
  		\left\{\begin{array}{l}
  			f^n_t+v^n\cdot\nabla f^n=g,\ x\in\mathbb{R}^d,\ t>0, \\
  			f^n(0,x)=f_0(x)
  		\end{array}\right.
  	\end{equation}
  	with ~$f_0\in B^s_{p,r},\ g\in L^1([0,T];B^s_{p,r})$, then there exists ~$\alpha\in L^1([0,T])$,  such that
  	$$\sup\limits_{n\in\bar{\mathbb{N}}}\|v^n(t)\|_{B^{s+1}_{p,r}}\leq \alpha(t).$$
  	If ~$v^n$ converges to ~$v^{\infty}$ in ~$L^1([0,T];B^s_{p,r})$ , then ~$f^n$ will converge to ~$f^{\infty}$ in ~$C([0,T];B^s_{p,r})$ .
  \end{lemm}
  \begin{lemm}\cite{Luo1}
  	Let $1\leq p\leq \infty$ and $1<r\leq\infty $. For any $\epsilon>0$, there exists $\gamma_0\in H^\infty$, such that the following holds:
  	
  	1. $\|\gamma_0\|_{B_{p.r}^{\frac{1}{p}}}\leq \epsilon$
  	\\
  	
  	2. There is a unique solution $\gamma\in C([0,T);H^\infty)$ to the equation\eqref{e1} with a maximal lifespan T <$\epsilon$;
  	\\
  	
  	3.$\limsup_{t\rightarrow T^-} \|\gamma\|_{B_{p.r}^{\frac{1}{p}}}\geq \limsup_{t\rightarrow T^-} \|\gamma\|_{B_{\infty,\infty}^{0}}=\infty$ 
  \end{lemm}

  \begin{lemm}\label{assumption}\cite{Yei}
	The space $B_{\infty,1}^0$ is not a Banach algebra. If we choose 
	\begin{align}
		u_0=[cos2^{N+5}x(1+N^{-\frac{1}{10}}S_Nh)+R_L^1]N^{-\frac{1}{10}}\label{initial}
	\end{align}
	where $R_L^1=-(-1-\partial_{xx})^{-1}[cos2^{N+5}x(1+N^{-\frac{1}{10}}S_Nh)]$ and $h=1_{x\geq0}$,
	
	then we have
	$$\|u_0\|_{B_{\infty,1}^0}\leq CN^{-\frac{1}{10}},~~~\|u_0\|_{B_{\infty,\infty,1}^0}\leq CN^{\frac{19}{10}},~~\|u_0^2\|_{B_{\infty,1}^0}\geq CN^{\frac{3}{5}}$$
\end{lemm}
\begin{lemm}\label{aaaa}
	For any $\Gamma \in B^0_{\infty,1}\cap B^0_{\infty,\infty,1}$, we have
	\begin{align}
		\|\Gamma \partial_xG^{-1}M\|_{B^0_{\infty,1}}&\leq C(\frac{1}{2\lambda}\|\Gamma\|^2_{B^0_{\infty,1}}\|\Gamma\|_{B^0_{\infty,\infty,1}}+\|\Gamma\|_{B^0_{\infty,\infty,1}}\|\Gamma\|_{B^0_{\infty,1}})
		\label{ie1}\\
		\|\Gamma \partial_xG^{-1}M\|_{B^0_{\infty\infty,1}}&\leq C(\frac{1}{2\lambda}\|\Gamma\|^2_{B^0_{\infty,1}}\|\Gamma\|_{B^0_{\infty,\infty,1}}+\|\Gamma\|_{B^0_{\infty,\infty,1}}\|\Gamma\|_{B^0_{\infty,1}})\label{ie2}
		\\
		&\|\Gamma^2\|_{B^0_{\infty,1}}\leq C\|\Gamma\|_{B^0_{\infty,1}}\|\Gamma\|_{B^0_{\infty,\infty,1}}\label{ie3}
		\\
		&\|\Gamma^2\|_{B^0_{\infty,\infty,1}}\leq C\|\Gamma\|_{B^0_{\infty,1}}\|\Gamma\|_{B^0_{\infty,\infty,1}}\label{ie4}
	\end{align}
\end{lemm}
\begin{lemm}\label{b}
	Define $R_{j}=G_{-1}M\Delta_j\Gamma_x-\Delta_j(G^{-1}M\Gamma_x), \tilde{R_{j}}=G^{-1}M\Delta_j(\frac{\Gamma}{2})_x-\Delta_j(G^{-1}M(\frac{\Gamma}{2})_x)$. Then we have 
	\begin{align}
		\sup_j ((j+2)^2\|R_{j}\|_{L^\infty})\leq C\|\Gamma\|_{B^0_{\infty,1}}(\|\Gamma\|_{B^0_{\infty,\infty,1}}+\frac{1}{2\lambda}\|\Gamma\|_{B^0_{\infty,\infty,1}}\|\Gamma\|_{B^0_{\infty,1}})
		\\
		\sum_j \|R_{j}\|_{L^\infty}\leq C\|\Gamma\|_{B^0_{\infty,1}}(\|\Gamma\|_{B^0_{\infty,\infty,1}}+\frac{1}{2\lambda}\|\Gamma\|_{B^0_{\infty,\infty,1}}\|\Gamma\|_{B^0_{\infty,1}})
		\\
		\sum_j \|\tilde{R_{j}}\|_{L^\infty}\leq C\big(\|\Gamma\|_{B^0_{\infty,1}}^2\|\Gamma\|_{B^0_{\infty,\infty,1}}+\frac{1}{2\lambda}(\|\Gamma\|_{B^0_{\infty,\infty,1}}\|\Gamma\|_{B^0_{\infty,1}})^2\big)
	\end{align}
\end{lemm}
\section{Some results for $B^0_{\infty,\infty,1}$}
~~~~~In this section, we give some detail of the calculations we will use in the Theorem \ref{the1}.

\textbf{Proof of Lemma \ref{aaaa}}
 To calculate \eqref{ie1} and \eqref{ie2}, we can rewrite $\Gamma \partial_xG^{-1}M$ as the form of $$\partial_xG^{-1}(\frac{1}{2\lambda}\Gamma^2+\Gamma_x)\Gamma,~~ G^{-1}M=G^{-1}(\frac{1}{2\lambda}\Gamma^2+\Gamma_x).$$ Using the fact that $\partial_x^2 G^{-1} \Gamma=G^{-1}\Gamma+\Gamma$, we have $\Gamma \partial_xG^{-1}M=G^{-1}(\frac{1}{2\lambda}\Gamma^2+\Gamma)+\Gamma$.
 Applying the Bony's decompositon to $\Gamma \partial_xG^{-1}M$, we can conclude that
 \begin{align}
 	\|\Gamma \partial_xG^{-1}M\|_{B_{\infty,1}^0}&\leq \|T_{\partial_xG^{-1}\frac{1}{2\lambda}\Gamma^2}\Gamma\|_{B_{\infty,1}^0}+\|T_{\Gamma}\partial_xG^{-1}\frac{1}{2\lambda}\Gamma^2\|_{B_{\infty,1}^0}+\|R(\Gamma,\partial_xG^{-1}\frac{1}{2\lambda}\Gamma^2)\|_{B_{\infty,1}^0}
 	\notag\\
 	&+2(2\|T_{\Gamma}\Gamma\|_{B_{\infty,1}^0}+\|R(\Gamma,\Gamma)\|_{B_{\infty,1}^0})
 	,
 \end{align} 
 and
 \begin{align}
 \| T_{\partial_xG^{-1}\Gamma^2}\Gamma\|_{B^0_{\infty,1}}=&\|\|\Delta_j(\sum S_{j^\prime-1}\partial_xG^{-1}\Gamma^2\Delta_{j^\prime}\Gamma)\|_{L^\infty}\|_{l^1}
 \\ \notag
 =&\sum_j\|\sum_{|j-j^\prime|\leq 4}\Delta_j(S_{j^\prime-1}\partial_xG^{-1}\Gamma^2\Delta_{j^\prime}\Gamma)\|_{L^\infty}
  \\ \notag
 \leq&\sum_j\sum_{|j-j^\prime|\leq 4}\|\Delta_j(S_{j^\prime-1}\partial_xG^{-1}\Gamma^2\Delta_{j^\prime}\Gamma)\|_{L^\infty}
 \\ \notag
 \leq&C\|S_{j^\prime-1}\partial_xG^{-1}\Gamma^2\|_{L^\infty}\|\Gamma\|_{B_{\infty,1}^0}.
 \end{align}
 For $\|S_{j^\prime-1}\partial_xG^{-1}\Gamma^2\|_{L^\infty}$, it is not hard to check that
 \begin{align}
 	\|S_{j^\prime-1}\partial_xG^{-1}\Gamma^2\|_{L^\infty}=&
 	\|\sum_{j^{\prime\prime}\leq j^\prime-2}\Delta_{j^{\prime\prime}}\partial_xG^{-1}\Gamma^2\|_{L^\infty}
 	\leq
 	 \sum_{j^{\prime\prime}\leq j^\prime-2}\frac{(j^{\prime\prime}+2)^2}{(j^{\prime\prime}+2)^2} \|\Delta_{j^{\prime\prime}}\partial_xG^{-1}\Gamma^2\|_{L^\infty}
 	 \notag \\
 	 &\leq C\|\partial_x G^{-1}\Gamma^2\|_{B_{\infty,\infty,1}^0}
 	 \leq C \|\Gamma^2\|_{B_{\infty,\infty,1}^0}.
 \end{align}
 Then we will have
 \begin{align}
 	\| T_{\partial_xG^{-1}\Gamma^2}\Gamma\|_{B^0_{\infty,1}}\leq C\|\Gamma^2\|_{B_{\infty,\infty,1}^0}\|\Gamma\|_{B_{\infty,1}^0}.\label{m1}
 \end{align}
By the same token, one has
 \begin{align}
 	\| T_\Gamma{\partial_xG^{-1}\Gamma^2}\|_{B^0_{\infty,1}}&=\|\|\Delta_j(\sum S_{j^\prime-1}\Gamma\Delta_{j^\prime}\partial_xG^{-1}\Gamma^2)\|_{L^\infty}\|_{l^1}
 	\notag \\
 	&\leq C\|\Gamma\|_{B_{\infty,\infty,1}^0}\|\Gamma^2\|_{B_{\infty,1}^0}.\label{m2}
 \end{align}
For remainder operator $R(\Gamma,\partial_xG^{-1}\Gamma^2)$, we have 
\begin{align}
	\|R(\Gamma,G^{-1}\Gamma^2)\|_{B_{\infty,1}^0}=&\|\|\Delta_{j^\prime}R(\Gamma,G^{-1}\Gamma^2)\|_{L^\infty}\|_{l^1}=\sum_{j^\prime}\|\Delta_{j^\prime}R(\Gamma,G^{-1}\Gamma^2)\|_{L^\infty}=\sum_{j^\prime}\|\Delta_{j^\prime} \sum_j R_j\|_{L^\infty}
	\notag\\
	&\leq C \sum_{j^\prime} \sum_{\substack{j\geq j^\prime-N_0 \\ |v|\leq 1}}\|\Delta_{j-v}\Gamma\Delta_j\partial_xG^{-1}\Gamma^2\|_{L^\infty}
	\notag \\
	&\leq C \sum_{j^\prime} \sum_{\substack{j \\ |v|\leq 1}}1_{j\geq j^\prime-N_0}\|\Delta_{j-v}\Gamma\|_{L^\infty}\|\Delta_j\partial_xG^{-1}\Gamma^2\|_{L^\infty}
	\notag \\
	&= C \sum_{j} \sum_{\substack{j^\prime \\ |v|\leq 1}}1_{j\geq j^\prime-N_0}\|\Delta_{j-v}\Gamma\|_{L^\infty}\|\Delta_j\partial_xG^{-1}\Gamma^2\|_{L^\infty}
	\notag\\
	&\leq C \sum_{\substack{j \\ |v|\leq 1}} (j+N_0+2)\|\Delta_{j-v}\Gamma\|_{L^\infty}\|\Delta_j\partial_xG^{-1}\Gamma^2\|_{L^\infty}
	\notag \\
	&\leq C \|\Gamma^2\|_{B_{\infty,\infty,1}^0}\|\Gamma\|_{B_{\infty,1}^0}.\label{m3}
\end{align}
For the estimates of $\|T_{\Gamma}\Gamma\|_{B_{\infty,1}^0},\|R(\Gamma,\Gamma)\|_{B_{\infty,1}^0}$, one can deduce that
\begin{align}
	\|T_{\Gamma}\Gamma\|_{B_{\infty,1}^0}&\leq C\|\Gamma\|_{B_{\infty,1}^0}\|\Gamma\|_{B_{\infty,\infty,1}^0}
	\notag \\
		\|R(\Gamma,\Gamma)\|_{B_{\infty,1}^0}&\leq C\|\Gamma\|_{B_{\infty,1}^0}\|\Gamma\|_{B_{\infty,\infty,1}^0}.
\end{align}
Then from the above estimates, we arrive at \eqref{ie3}. Similarly, we can prove \eqref{ie4} and here we omit it. Utilising \eqref{ie3},\eqref{ie4} and \eqref{m1}-\eqref{m3}, we conclude that
\begin{align}
	\| T_{\partial_xG^{-1}\Gamma^2}\Gamma\|_{B^0_{\infty,1}}&\leq C\|\Gamma\|_{B_{\infty,\infty,1}^0}\|\Gamma\|_{B_{\infty,1}^0}^2
	\notag\\
		\| T_\Gamma{\partial_xG^{-1}\Gamma^2}\|_{B^0_{\infty,1}}&\leq C\|\Gamma\|_{B_{\infty,\infty,1}^0}\|\Gamma\|_{B_{\infty,1}^0}^2
		\notag\\
		\|R(\Gamma,G^{-1}\Gamma^2)\|_{B_{\infty,1}^0}&\leq C\|\Gamma\|_{B_{\infty,\infty,1}^0}\|\Gamma\|_{B_{\infty,1}^0}^2
\end{align}
Thus we complete the proof of \eqref{ie1}. The proof of \eqref{ie2} is very similar, and we shall omit it.

$\hfill\square$

\textbf{Proof of Lemma \ref{b}}
 We will prove Lemma\ref{b} with similar calculations as Lemma 100 in \cite{Chemin2011} and here we present some details of it.
 Our main difficulty is to estimate  $R_j = G^{-1}M\Delta_j \Gamma_x-\Delta_j(G^{-1}M\Gamma_x)$. Applying Bony's decomposition, we have $R_j=\sum_{j=1}^{8}$, where 
 \begin{align}
 	R_j^1&=[T_{{G^{-1}\tilde M}},\Delta_j]\partial_x\Gamma, ~~~~R_j^2=T_{\partial_x\Delta_j\Gamma} {G^{-1}\tilde M},
 	\notag\\
 	R_j^3&=-\Delta_jT_{\partial_x\Gamma} {G^{-1}\tilde M}, ~~~~R_j^4=\partial_x R(G^{-1}\tilde{M},\Delta_j\Gamma),
 	\notag\\
 	R_j^5&=-R(\partial_xG^{-1}\tilde{M},\Delta_j\Gamma), ~~~~R_j^6=\partial_x \Delta_jR(G^{-1}\tilde{M},\Gamma),
 	\notag\\
 	R_j^7&=-\Delta_jR(G^{-1}\tilde{M},\Gamma), ~~~~R_j^8=[S_0G^{-1}M,\Delta_j]\partial_x\Gamma,
 	\notag\\
 	\tilde{M}&=M-S_0M.
 \end{align}
Firstly, by Lemma \ref{prop}, we can deduce that
\begin{align}
	R_j^1=[T_{{G^{-1}\tilde M}},\Delta_j]\partial_x\Gamma=\sum_{|j-j^\prime|\leq 4}[S_{j^\prime-1}G^{-1}\tilde{M},\Delta_j]\Delta_{j^\prime}\Gamma_x,
\end{align}   
hence 
\begin{align}
	\sum_{j}\|R_j^1\|_{L^\infty} &\leq C\sum_{j}\sum_{|j-j^\prime|\leq 4}\|\partial_x S_{j^\prime-1}G^{-1}\tilde{M}\|_{L^\infty}\|\Delta_{j^\prime}\Gamma\|_{L^\infty}
	\notag\\
	&\leq C\|\partial_xG^{-1}M\|_{B_{\infty,\infty,1}^0} \|\Gamma\|_{B_{\infty,1}^0} 
	\notag\\
	&\leq C(\|\Gamma\|_{B_{\infty,\infty,1}^0}+\frac{1}{2\lambda}\|\Gamma\|_{B_{\infty,1}^0}\|\Gamma\|_{B_{\infty,\infty,1}^0})\|\Gamma\|_{B_{\infty,1}^0} .
\end{align}
Consequently, one can get
\begin{align}
	\|(j+2)^2\|R_j^1\|_{L^\infty}\|_{l^\infty} &\leq C\|(j+2)^2\sum_{|j^\prime-j|\leq 4}\|\partial_xS_{j^\prime-1}G^{-1}\tilde{M}\|_{L^\infty}\|\Delta_j\Gamma\|_{L^\infty}\|_{l^\infty}
	\notag\\
	&\leq C \|(j+2)^2\sum_{|j-j^\prime|\leq 4} \|\partial_xS_{j^\prime-1}G^{-1}\frac{\Gamma^2}{2\lambda}\|_{L^\infty}\|\Delta_{j'}\Gamma\|_{L^\infty}\|_{l^\infty}
	\notag\\
	&+C \|(j+2)^2\sum_{|j-j^\prime|\leq 4} \|S_{j^\prime-1}\Gamma\|_{L^\infty}\|\Delta_{j'}\Gamma\|_{L^\infty}\|_{l^\infty}
	\notag\\
	&\leq C \|(j+2)^2 \sum_{|j-j^\prime|\leq 4}\sum_{j^{\prime\prime}\leq j^\prime-2}\frac{1}{2^{j^\prime}}\|\partial_x^2G^{-1}\Delta_{j^{\prime\prime}}\frac{\Gamma^2}{2\lambda}\|_{L^\infty}\|_{l^\infty}
	+C\|\Gamma\|_{B_{\infty,\infty,1}^0}\|\Gamma\|_{B_{\infty,1}^0}
	\notag\\
	&\leq C(\| \sum_{|j-j^\prime|\leq 4}\sum_{j^{\prime\prime}\leq j^\prime-2}\frac{(j+2)^2}{2^{j^\prime}}\|_{l^\infty}\frac{\|\Gamma^2\|_{B_{\infty,1}^0}\|\Gamma\|_{B_{\infty,1}^0}}{2\lambda}+\|\Gamma\|_{B_{\infty,\infty,1}^0}\|\Gamma\|_{B_{\infty,1}^0})
	\notag\\
	&\leq C(\frac{\|\Gamma\|_{B_{\infty,1}^0}^2\|\Gamma\|_{B_{\infty,\infty,1}^0}}{2\lambda}+\|\Gamma\|_{B_{\infty,\infty,1}^0}\|\Gamma\|_{B_{\infty,1}^0}).
\end{align}
For the estimate of $R_j^2$, we have 
$$R_j^2=T_{\partial_x\Delta_j\Gamma}G^{-1}\tilde{M},$$
then
\begin{align}
	\|\|R_j^2\|_{L^\infty}\|_{l^1}&\leq C\|\sum_{j^{\prime}\geq j-3}\|S_{j^\prime-1}\partial_x\Delta_j\Gamma\Delta_{j^\prime}G^{-1}\tilde{M}\|_{L^\infty}\|_{l^1}
	\notag\\
	&\leq C\|\sum_{j^{\prime}\geq j-3}\|S_{j^\prime-1}\Delta_j\Gamma\|_{L^\infty}\|\Delta_{j^\prime}G^{-1}\partial_x\frac{\Gamma^2}{2\lambda}\|_{L^\infty}\|_{l^1}+\|\sum_{j^{\prime}\geq j-3}\|S_{j^\prime-1}\partial\Delta_j\Gamma\|_{L^\infty}\|\Delta_{j^\prime}G^{-1}\Gamma\|_{L^\infty}\|_{l^1}
	\notag\\
	&\leq C (\|\sum_{j^{\prime}\geq j-3}\|S_{j^\prime-1}\Delta_j\Gamma\|_{L^\infty}\|\Delta_{j^\prime}G^{-1}\partial_x\frac{\Gamma^2}{2\lambda}\|_{L^\infty}\|_{l^1}+\|\Gamma\|_{B_{\infty,1}^0}\|\Gamma\|_{B_{\infty,\infty,1}^0}
	)
	\notag\\
	&\leq C \|\Gamma\|_{B_{\infty,1}^0}(\|\sum_{j^{\prime}\geq j-3}\frac{1}{2^{j^\prime}}(\|\Delta_{j^\prime}\frac{\Gamma^2}{2\lambda}\|_{L^\infty}+\|\Delta_{j^\prime}G^{-1}\frac{\Gamma^2}{2\lambda}\|_{L^\infty})\|_{l^1}+\|\Gamma\|_{B_{\infty,\infty,1}^0})
	\notag\\
	&\leq C \|\Gamma\|_{B_{\infty,1}^0}(\sum_j \sum_{j^{\prime}\geq j-3}\frac{1}{2^{j^\prime}}(\|\Delta_{j^\prime}\frac{\Gamma^2}{2\lambda}\|_{L^\infty}+\|\Delta_{j^\prime}G^{-1}\frac{\Gamma^2}{2\lambda}\|_{L^\infty})+\|\Gamma\|_{B_{\infty,\infty,1}^0})
	\notag\\
	&\leq C \|\Gamma\|_{B_{\infty,1}^0}(\sum_j \sum_{j^{\prime}\geq j-3}\frac{1}{2^{j^\prime}}c_{j^\prime}\|\frac{\Gamma^2}{2\lambda}\|_{B_{\infty,1}^0}+\|\Gamma\|_{B_{\infty,\infty,1}^0})
	\notag \\
	&\leq C \|\Gamma\|_{B_{\infty,1}^0}(\frac{1}{2\lambda}\|\Gamma\|_{B_{\infty,1}^0}\|\Gamma\|_{B_{\infty,\infty,1}^0}+\|\Gamma\|_{B_{\infty,\infty,1}^0}) ,
\end{align}
where $\{c_{j^\prime}\}_{j^\prime\geq-1}$ is a sequence such that $\|c_{j^\prime}\|_{l^1}\leq 1$.
In addition, we have
\begin{align}
	\|(j+2)^2\|R_j^2\|_{L^\infty}\|_{l^\infty}=&\|(j+2)^2\|\sum_{j^{\prime}\geq j-3}S_{j^\prime-1}\partial_x\Delta_j\Gamma\Delta_{j^\prime}G^{-1}\tilde{M}\|_{L^\infty}\|_{l^\infty}
	\notag\\
	&\leq C (\|\Gamma\|_{L^\infty} \|(j+2)^2\sum_{j^\prime\geq j-3}\|\Delta_{j^\prime}\partial_xG^{-1}\frac{\Gamma^2}{2\lambda}\|_{L^\infty}\|_{l^\infty}+\|\Gamma\|_{B_{\infty,1}}\|\Gamma\|_{B_{\infty,\infty,1}})
	\notag\\
	&\leq C\|\Gamma\|_{B_{\infty,1}^0}(\|(j+2)^2\sum_{j^{\prime}\geq j-3}2^{-j^\prime}\|\Delta_{j^\prime}\frac{\Gamma^2}{2\lambda}+\Delta_{j^\prime}G^{-1}\frac{\Gamma^2}{2\lambda}\|_{L^\infty}\|_{l^\infty}+\|\Gamma\|_{B_{\infty,\infty,1}^0})
	\notag\\
	&\leq C\|\Gamma\|_{B_{\infty,1}^0}(\|\frac{\Gamma^2}{2\lambda}\|_{B_{\infty,\infty,1}^0}+\|\Gamma\|_{B_{\infty,1}^0}).
\end{align}
By Propsition\ref{prop}, we obatin the following estimate,
\begin{align}
	\sum_{j}\|R_j^3\|_{L^\infty}&\leq C \sum_j\sum_{\substack{|j-j^\prime|\leq4 \\ j^{\prime\prime}\leq j^\prime-1}}\|\Delta_j(\Delta_{j^{\prime\prime}}\Gamma_x\Delta_{j^\prime}G^{-1}M)\|_{L^\infty}
	\notag\\
	&\leq C \sum_j\sum_{\substack{|j-j^\prime|\leq4 \\ j^{\prime\prime}\leq j^\prime-1}}\|\Delta_{j^{\prime\prime}}\Gamma_x\|_{L^\infty}\|\Delta_{j^\prime}G^{-1}M\|_{L^\infty}
	\leq C \sum_j\sum_{\substack{|j-j^\prime|\leq4 \\ j^{\prime\prime}\leq j^\prime-1}}\|\Delta_{j^{\prime\prime}}\Gamma\|_{L^\infty}\|\partial_x\Delta_{j^\prime}G^{-1}M\|_{L^\infty}
	\notag \\
	&\leq C \sum_j\sum_{\substack{|j-j^\prime|\leq4 \\ j^{\prime\prime}\leq j^\prime-1}}(\|\Delta_{j^{\prime\prime}}\Gamma\|_{L^\infty}\|\Delta_{j^\prime}\partial_xG^{-1}\frac{\Gamma^2}{2\lambda}\|_{L^\infty}+\|\Delta_{j^{\prime\prime}}\Gamma\|_{L^\infty}\|\Delta_{j^{\prime}}\Gamma\|_{L^\infty})
	\notag\\
	&\leq C \|\Gamma\|_{B_{\infty,1}^0}(\sum_{j}\sum_{|j-j^\prime|\leq 4}\|\partial_x\Delta_{j^\prime}G^{-1}\frac{\Gamma^2}{2\lambda}\|_{L^\infty}+\|\Gamma\|_{B_{\infty,\infty,1}^0})
	\notag\\
	&\leq C \|\Gamma\|_{B_{\infty,1}^0}(\frac{\|\Gamma\|_{B_{\infty,\infty,1}^0}\|\Gamma\|_{B_{\infty,1}^0}}{2\lambda}+\|\Gamma\|_{B_{\infty,\infty,1}^0}).
\end{align}
Similarly, one arrive at
\begin{align}
	\|(j+2)^2\|R_j^3|_{L^\infty}\|_{l^\infty}&=\|(j+2)^2\|\sum_{\substack{|j-j^\prime|\leq 4 \\ j^{\prime\prime}\leq j^\prime-2}}\Delta_j(\Delta_{j^{\prime\prime}}\partial_x\Gamma\Delta_{j^{\prime}})G^{-1}M\|_{L^\infty}\|_{l^\infty}
		\notag\\
		&\leq C\|(j+2)^2\sum_{\substack{|j-j^\prime|\leq 4 \\ j^{\prime\prime}\leq j^\prime-2}}\|\Delta_{j^{\prime\prime}}\partial_x\Gamma\Delta_{j^\prime}G^{-1}M\|_{L^\infty}\|_{l^\infty}
		\notag\\
		&\leq C\|\Gamma\|_{B_{\infty,1}^0}\|G^{-1}M\|_{B_{\infty,\infty,1}^0}
		\notag\\
			&\leq C \|\Gamma\|_{B_{\infty,1}^0}(\frac{\|\Gamma\|_{B_{\infty,\infty,1}^0}\|\Gamma\|_{B_{\infty,1}^0}}{2\lambda}+\|\Gamma\|_{B_{\infty,\infty,1}^0}),
	\end{align}
The rest of the proof is very similar to the proof of Lemma100 in \cite{Chemin2011} and here we omit it.
$\hfill\square$
\begin{rema}
	We notice that in the proof of above two lemmas we utilised the special structure of the $G^{-1}M$, which contains a quadratic term of $\Gamma$. The quardratic term will generate a quadratic term of $B_{\infty,\infty,1}^0$, which will bring us with massive difficulties to bound the norm of $\|\Gamma(T)\|_{L^\infty_t(B_{\infty,\infty,1}^0)}$, which we will see later. Taking advantage of Bernstein's inequality \eqref{Bernstein} to transform the function in a new form, then constructing a convergent sequence will provide us with a easy way to get the estimate.
\end{rema}
\section{Local well-posedness}
~~~~In this section, we mainly study the local well-posedness for the modified Camassa-Holm equation. To prove Theorem\ref{them2}, we must recall a useful lemma. Firstly, we caonsider the following Cauchy problem for a general abstract equation 
\begin{equation}\label{eq10001}
	\left\{\begin{array}{l}
		u_t+A(u)u_x=F(u), t>0,\quad  x\in\mathbb{R} \\
		u(0,x)=u_0(x),  \quad x\in\mathbb{R}.
	\end{array}\right.
\end{equation}

Where A(u) is a polynomial of u, F is called a 'good operator' such that for any $\phi\in\mathcal{C}^\infty_0$ and any $\epsilon>0$ small enough,
\begin{align}
	if\quad u_n\phi\rightarrow u\phi\quad in \quad B^{1+\frac{1}{p}-\epsilon}_{p,1}
\end{align}
then
\begin{align}
	\left\langle F(u_n,\phi)\right\rangle\rightarrow \left\langle F(u),\phi \right\rangle
\end{align}
 The associated Lagrangian scale of \eqref{eq10001} is the following initial value problem
 \begin{equation}\label{liu}
 	\left\{\begin{array}{l}
 		y_t(t,\xi)=A(u)(t,y)\quad ,t\in [0,T)  \\
 		
 		y(0,\xi)=\xi,\quad \xi\in\mathbb{R}.
 	\end{array}\right.
 \end{equation}
 Introducing the new variable $U(t,\xi)=u(t,y(t,\xi))$. Then, combining \eqref{eq10001} and \eqref{liu}, we deduce that $u(t,\xi)$ satisfies the following equation.
 \begin{equation}
 	\left\{\begin{array}{l}
 		U_t(t,\xi)=F(u)(t,y(t,\xi))\triangleq \tilde{F}(U,y) \quad ,t\in [0,T)  \\
 		
 		U(0,\xi)=u_0(\xi),\quad \xi\in\mathbb{R}.
 	\end{array}\right.
 \end{equation}
Here is the Lemma:
\begin{lemm}\cite{2021}
	Let $u_0\in B^{1+\frac{1}{p}}_{p,1}$ with 1$\leq$p<$\infty$ and k$\in \mathbb{N}^+$. Suppose F is a 'good operator' and F,$\tilde{F}$ satisfy the following conditions:
	\begin{align}\label{F1}
		\|F(u)\|_{B_{p,1}^{\frac{1}{p}+1}}&\leq C(\|u\|_{B_{p,1}^{\frac{1}{p}+1}}^{k+1}+1),
    \end{align}
\begin{align}\label{F2}
		\|\tilde{F}(U,y)-\tilde{F}(\bar{U},\bar{y})\|_{W^{1,\infty}\cap W^{1,\infty}}&\leq C(\|U-\bar{U}\|_{W^{1,\infty}\cap W^{1,\infty}}+\|y-\bar{y}\|_{W^{1,\infty}\cap W^{1,\infty}}) ,
\end{align}
\begin{align}
		\|F(u)-F(\bar{u})\|_{B_{p,1}^{\frac{1}{p}+1}}&\leq C\|u-\bar{u}\|_{B_{p,1}^{\frac{1}{p}+1}}(\|u\|_{B_{p,1}^{\frac{1}{p}+1}}^{k}+\|\bar{u}\|_{B_{p,1}^{\frac{1}{p}+1}}^{k}+1),\label{F3}
	\end{align}
Then , there exists a time T \textgreater0 such that

(1) Existence: If \eqref{F1} holds, then \eqref{eq10001} has a solution u$\in E^P_T\triangleq \mathcal{C}_T(B_{p,1}^{\frac{1}{p}+1}\cap \mathcal{C}^1_T(B_{p,1}^{\frac{1}{p}})$;

(2)Uniqueness: If \eqref{F1} and \eqref{F2}, then the solution of \eqref{eq10001} is unique;

(3)Continuous dependence: If \eqref{F1}-\eqref{F3} hold, then the solution map is continuous from any bounded subset of $B_{p,1}^{\frac{1}{p}+1}$ to ${C}_T(B_{p,1}^{\frac{1}{p}+1})$

That is, the problem \eqref{eq10001} is locally well-posed in the sense of Hadamard. 
\end{lemm}

\textbf{Proof of Theorem \eqref{them2}}
For existence, since $\gamma^0 \in B^{0}_{\infty,1}(\mathbb{R})\cap B^{0}_{\infty,\infty,1}(\mathbb{R})$, we can deduce that 
\begin{align}
	\|G^{-1}m^n\|_{B^{0}_{\infty,1}\cap B^{0}_{\infty,\infty,1}}&\leq C\|m^n\|_{B^{-1}_{\infty,1}\cap B^{-1}_{\infty,\infty,1}}
	\notag  \\
	&\leq C(\|\gamma^n\|_{B^{0}_{\infty,1}\cap B^{0}_{\infty,\infty,1}}+
	\|\gamma^n\|^2_{B^{0}_{\infty,1}\cap B^{0}_{\infty,\infty,1}}), \label{2ineq1}
\end{align}
and
\begin{align}
	&\|\frac {(\gamma^{n})^2}{2}\ +\lambda G^{-1}m^{n}-\gamma G^{-1}m^{n}_x\|_{B^{0}_{\infty,1}\cap B^{0}_{\infty,\infty,1}}\\
	&\leq  C(\|\gamma^n\|^2_{B^{0}_{\infty,1}\cap B^{0}_{\infty,\infty,1}}+\|G^{-1}m^n\|_{B^{0}_{\infty,1}\cap B^{0}_{\infty,\infty,1}}
	+\|\gamma^n\|_{B^s_{p,r}}\|G^{-1}m_x^n\|_{B^{0}_{\infty,1}\cap B^{0}_{\infty,\infty,1}}) \notag  \\
	&\leq C(\|\gamma^n\|_{B^{0}_{\infty,1}\cap B^{0}_{\infty,\infty,1}}
	+\|\gamma^n\|^2_{B^{0}_{\infty,1}\cap B^{0}_{\infty,\infty,1}}
	+\|\gamma^n\|^3_{B^{0}_{\infty,1}\cap B^{0}_{\infty,\infty,1}}). \label{2ineq2}
\end{align}

Then we use the compactness method can easily get the result. For uniqueness, with the good structure of the transport term, we lift regularity by using operator $\partial_x-1$. We will transform\eqref{eq1} to a better form about n to obtain the uniqueness of $n$, then the uniqueness of $\gamma$ can be deduced from the uniqueness of $n$. For
the continuous dependence, the first two steps provide us with great convenience to the problem. The proof is very similar to the proof in \cite{He} and we omit it.
\section{Ill-posedness}

~~~~~In this section, we study the ill-posedness for the Cauchy problem\eqref{eq1} of the modified Camassa-Holm in $B_{\infty,1}^0$.

\begin{proof}\
	Let $\gamma$ be a solution to the modified Camassa-Holm equation with the initial data $\gamma_0$ defined as \eqref{initial}. Set 
 \begin{equation}\label{flow}
	\left\{\begin{array}{l}
	\frac{d}{dt} y(t,\xi)=G^{-1}m(t,y(t,\xi))  ,t\in [0,T)  \\
		
	y_0(\xi)=\xi,\quad \xi\in\mathbb{R}.
	\end{array}\right.
\end{equation}

   Therefore according to \eqref{flow}, we can find a time $T_0 > 0$ sufficiently small such that $\frac{1}{2}\leq y_\xi(t) \leq2$, for any t$\in[0,T_0]$. Let $T=N^{-\frac{1}{2}}\leq T_0$ for $N > 10$ large enough. To prove the norm inflation, it suffices to prove there exists a time $t_0\in[0,N^{-\frac{1}{2}}]$ such that $\|\gamma(t_0)\|_{B_{\infty,1}^0}\geq \ln N$ for $N>10$ large enough. Let us assume the opposite. Namely, we suppose that
   \begin{align}\label{as}
   	\sup_{t\in[0,N^{-\frac{1}{2}}]}\|\gamma(t)\|<\ln N
   \end{align}
 Applying $\Delta_j$ to\eqref{eq1}
 \begin{align}
 	\Delta_j \gamma_t+G^{-1}m\Delta_j\gamma_x=\Delta_j\frac{\gamma^2}{2}+\Delta_j\lambda G^{-1}m-\Delta_j(\gamma G^{-1}m_x)+G^{-1}m\Delta_j\gamma_x-\Delta_j(G^{-1}m\gamma_x),
 \end{align}
 and we denote $\Gamma(t,\xi)\triangleq\gamma\circ y,M(t,\xi)\triangleq m\circ y$, we get
 \begin{align}\label{01}
 	\Delta_j\Gamma_t=\Delta_j (\frac{\Gamma^2}{2})+\Delta_j \lambda G^{-1}M-\Delta_j(\Gamma G^{-1}M_x)+R_j\circ y,
 \end{align}
 with $$R_j\circ y \triangleq G^{-1}M\Delta_j \Gamma_x-\Delta_j(G^{-1}M\Gamma_x).$$
 Integrating \eqref{01} with respect to t, we obtain
 \begin{align}
 	\Delta_j\Gamma(t)&-\Delta_j\Gamma(0)=\int_{0}^{t}\Delta_j (\frac{\Gamma^2}{2})(t^{\prime})+\Delta_j (\lambda G^{-1})M(t^{\prime})-\Delta_j(\Gamma G^{-1}M_x)(t^{\prime})+R_j\circ y(t^{\prime}) dt^{\prime}
 	\notag \\
 	\Delta_j\Gamma(t)&=\Delta_j\Gamma(0)+\int_{0}^{t}\Delta_j (\frac{\Gamma^2}{2})(t^{\prime})+\Delta_j (\lambda G^{-1})M(t^{\prime})-\Delta_j(\Gamma G^{-1}M_x)(t^{\prime})+R_j\circ y(t^{\prime}) dt^{\prime}
 	\notag \\
 	&=\Delta_j\Gamma(0)+\int_{0}^{t}\Delta_j (\frac{\Gamma^2}{2})(0)+\Delta_j (\frac{\Gamma^2}{2})(t^{\prime})-\Delta_j (\frac{\Gamma^2}{2})(0)+\Delta_j (\lambda G^{-1})M(t^{\prime})
 	\notag \\
 	&-\Delta_j(\Gamma G^{-1}M_x)(t^{\prime})+R_j\circ y(t^{\prime}) dt^{\prime},
 \end{align}
According to the above inequalities, one has
\begin{align}\label{keykey}
	\sup \limits_{t\in [0,T]}\|\Gamma(t)\|_{B_{\infty,1}^0}& \geq t\|\frac{\Gamma_0^2}{2}\|_{B_{\infty,1}^0}-\|\Gamma_0\|_{B_{\infty,1}^0}-\int_{0}^{t}\|\frac{\Gamma^2}{2}-\frac{\Gamma_0^2}{2}\|_{B_{\infty,1}^0}dt^{\prime}-\int_{0}^{t}\lambda\|G^{-1}M\|_{B_{\infty,1}^0}dt^{\prime}
	\notag\\
	&-\int_{0}^{t}\|\Gamma G^{-1}M_x\|_{B_{\infty,1}^0}dt^{\prime}-\int_{0}^{t}
	\sum_{j}\|R_j\circ y\|_{L^\infty}dt^{\prime}
\end{align}
Multiple \eqref{01} with $\gamma$, we have 
\begin{align}\label{03}
	\frac{1}{2}\frac{d}{dt} {\gamma}^2+\frac{1}{2}G^{-1}m(\gamma)_x^2=\frac{{\gamma}^3}{2}+\lambda \gamma G^{-1}m-\gamma^2G^{-1}m_x,
\end{align}
Applying $\Delta_j$ and Lagrange coordinates to \eqref{03} yields
\begin{align}\label{04}
	\frac{1}{2}\frac{d}{dt} \Delta_j \Gamma^2=\Delta_j\frac{{\Gamma}^3}{2}+\lambda\Delta_j(\Gamma G^{-1}M)-\Delta_j(\Gamma^2 G^{-1}M_x)+\tilde{R_j}
\end{align}
with $\tilde{R_j} \triangleq G^{-1}M\Delta_j(\Gamma_x\Gamma)-\Delta_j (G^{-1}M\Gamma_x\Gamma)$

Integrating \eqref{04} with respect to t, we infer that
\begin{align}
	\Delta_j(\frac{\Gamma^2}{2} (t)-\frac{\Gamma^2_0}{2})=\int_{0}^{t} \Delta_j\frac{{\Gamma}^3}{2}(\tau)+\lambda\Delta_j(\Gamma G^{-1}M)(\tau)-\Delta_j(\Gamma^2 G^{-1}M_x)(\tau)+\tilde{R_j}(\tau)d\tau,
\end{align}
from which we deduce that
\begin{align}
	\|\frac{\Gamma^2}{2} (t)-\frac{\Gamma^2_0}{2}\|_{B_{\infty,1}^0}\leq C\int_{0}^{t} \|\frac{{\Gamma}^3}{2}\|_{B_{\infty,1}^0}+\lambda\|(\Gamma G^{-1}M)\|_{B_{\infty,1}^0}+\|(\Gamma^2 G^{-1}M_x)\|_{B_{\infty,1}^0}+\|\tilde{R_j}\|_{B_{\infty,1}^0}dt.
\end{align}
Together with Lemma\ref{aaaa} and Lemma \ref{b}, we arrive at
\begin{align}\label{key1}
	\|\frac{\Gamma^2}{2} (t)-\frac{\Gamma^2_0}{2}\|_{B_{\infty,1}^0}&\leq C \int_{0}^{t} \|\Gamma(\tau)\|_{B_{\infty,\infty,1}^0}^2\|\Gamma(\tau)\|_{B_{\infty,1}^0}+\|\Gamma(\tau)\|_{B_{\infty,1}^0}\|\Gamma(\tau)\|_{B_{\infty,\infty,1}^0}^2
	\notag\\
	&+\|\Gamma(\tau)\|_{B^0_{\infty,\infty,1}}\|\Gamma(\tau)\|_{B^0_{\infty,1}}+\|\Gamma(\tau)\|_{B^0_{\infty,\infty,1}}^2\|\Gamma(\tau)\|_{B^0_{\infty,1}}^2 d\tau.
\end{align}
Following the similar proof of Lemma\ref{existence} and Lemma\ref{priori estimate}, we deduce
\begin{align}
	\|\Gamma(T)\|_{L^\infty_t(B_{\infty,\infty,1}^0)} &\leq \|\Gamma(T)\|_{L^\infty_t(B_{\infty,1}^0\cap B_{\infty,\infty,1}^0)}
	\notag\\
	&\leq \|\Gamma(0)\|_{B_{\infty,1}^0\cap B_{\infty,\infty,1}^0}+\int_{0}^{T}\|\frac{\Gamma^2(t)}{2}\|_{B_{\infty,1}^0\cap B_{\infty,\infty,1}^0}+\|\lambda G^{-1}M(t)\|_{B_{\infty,1}^0\cap B_{\infty,\infty,1}^0}
	\notag\\
	&~~+\|\Gamma G^{-1}M_x(t)\|_{B_{\infty,1}^0\cap B_{\infty,\infty,1}^0}+\sum_{j}\|R_j(t)\|_{L^\infty}+\sup_j((j+2)^2\|R_j(t)\|_{L^\infty}) dt
	\notag\\
	&\leq \|\Gamma(0)\|_{B_{\infty,1}^0\cap B_{\infty,\infty,1}^0}+C\int_{0}^{T}\frac{\|\Gamma(t)\|_{B_{\infty,1}^0}\|\Gamma(t)\|B_{\infty,\infty,1}^0}{2} +\lambda \|\Gamma\|_{B_{\infty,1}^0\cap B_{\infty,\infty,1}^0}
	\notag\\
	&
	~~+\lambda\|\Gamma(t)\|_{B_{\infty,1}^0}\|\Gamma(t)\|B_{\infty,\infty,1}^0+\frac{1}{2\lambda}\|\Gamma\|^2_{B^0_{\infty,1}}\|\Gamma\|_{B^0_{\infty,\infty,1}}+\|\Gamma\|_{B^0_{\infty,\infty,1}}\|\Gamma\|_{B^0_{\infty,1}}dt,
\end{align}
by our assumption and Lemma\eqref{assumption}
\begin{align}
	\|\Gamma(T)\|_{L^\infty_t(B_{\infty,\infty,1}^0)}&\leq CN^{\frac{19}{10}}
		+N^{-\frac{1}{2}}\|\Gamma(T)\|_{L^\infty_t(B_{\infty,\infty,1}^0)}+N^{-\frac{1}{2}}\ln N
	\notag\\
	&~~+N^{-\frac{1}{2}}\|\Gamma(T)\|_{L^\infty_t(B_{\infty,\infty,1}^0)}+\lambda N^{-\frac{1}{2}}\ln N+\frac{1}{2}N{-\frac{1}{2}}\ln^2N\|\Gamma(T)\|_{L^\infty_t(B_{\infty,\infty,1}^0)},
\end{align}
with N sufficiently large. 
  
  Consequently, we have
\begin{align}\label{key}
	\|\Gamma(T)\|_{L^\infty_t(B_{\infty,\infty,1}^0)}&\leq CN^{\frac{19}{10}}.
\end{align}
Taking advantage of \eqref{key} and our assumption, we can deduce that \eqref{key1} can be rewritten as 
\begin{align}
	\|\frac{\Gamma^2}{2} (t)-\frac{\Gamma^2_0}{2}\|_{B_{\infty,1}^0} \leq&~~C~(N^{-\frac{1}{2}}N^{\frac{38}{10}}\ln N+N^{-\frac{1}{2}}N^{\frac{19}{10}}\ln^2 N
		\notag\\
		&+N^{\frac{14}{10}}\ln N+N^{-\frac{1}{2}}N^{\frac{38}{10}}\ln^2 N)
\end{align}
Utilising Lemma\ref{aaaa} and Lemma\ref{b}, one has
\begin{align}
	\int_{0}^{t}\lambda\|G^{-1}M\|_{B_{\infty,1}^0}dt^{\prime}\leq C (N^{-\frac{1}{2}}\ln N+N^{-\frac{1}{2}}N^{\frac{19}{10}}\ln N), \label{key2}
\end{align}
\begin{align}
		\int_{0}^{t}\|\Gamma G^{-1}M_x\|_{B_{\infty,1}^0}dt^{\prime}\leq C N^{-\frac{1}{2}}(N^{-\frac{19}{10}}\ln^2 N +N^{\frac{19}{10}}\ln N),\label{key3}
\end{align}

\begin{align}
\int_{0}^{t}
	\sum_{j}\|R_j\circ y\|_{L^\infty}dt^{\prime} \leq CN^{-\frac{1}{2}}(N^{-\frac{19}{10}}\ln^2 N +N^{\frac{19}{10}}\ln N).\label{key4}
\end{align}

Inserting \eqref{key1} and \eqref{key2}-\eqref{key4} into \eqref{keykey}, we deduce that for any t$\in[0,T]$
\begin{align}
	\|\gamma(t)\|_{B_{\infty,1}^0}=\sum_{j}\|\Delta_j\gamma\|_{L^\infty}&=\sum_{j}\|\gamma\circ y\|_{L^\infty}
	\notag\\
	&\geq t\|\frac{\Gamma_0^2}{2}\|_{B_{\infty,1}^0}-N^{-\frac{1}{10}}-C(N^{-\frac{1}{2}}N^{\frac{38}{10}}\ln N+N^{-\frac{1}{2}}N^{\frac{19}{10}}\ln^2 N+N^{\frac{14}{10}}\ln N
	\notag\\
	&~~+N^{-\frac{1}{2}}N^{\frac{38}{10}}\ln^2 N)-C (N^{-\frac{1}{2}}\ln N+N^{-\frac{1}{2}}N^{\frac{19}{10}}\ln N)
	\notag\\
	&~~-C N^{-\frac{1}{2}}(N^{-\frac{19}{10}}\ln^2 N +N^{\frac{19}{10}}\ln N)-CN^{-\frac{1}{2}}(N^{-\frac{19}{10}}\ln^2 N +N^{\frac{19}{10}}\ln N)
	\notag\\
	&\geq CtN^{\frac{3}{5}}-N^{-\frac{1}{10}}-C(N^{-\frac{1}{2}}N^{\frac{38}{10}}\ln N+N^{-\frac{1}{2}}N^{\frac{19}{10}}\ln^2 N+N^{\frac{14}{10}}\ln N
	\notag\\
	&~~+N^{-\frac{1}{2}}N^{\frac{38}{10}}\ln^2 N)-C (N^{-\frac{1}{2}}\ln N+N^{-\frac{1}{2}}N^{\frac{19}{10}}\ln N)
	\notag\\
	&~~-C N^{-\frac{1}{2}}(N^{-\frac{19}{10}}\ln^2 N +N^{\frac{19}{10}}\ln N)-CN^{-\frac{1}{2}}(N^{-\frac{19}{10}}\ln^2 N +N^{\frac{19}{10}}\ln N)
	\notag\\
	&\geq CtN^{\frac{3}{5}},
\end{align}
where the second inequality holds by Lemma\eqref{assumption}. That is 
\begin{align}
	\|\gamma(t)\|_{B_{\infty,1}^0}\geq CN^{\frac{3}{5}-\frac{1}{2}}, \forall t\in[\frac{1}{2N^\frac{1}{2}},\frac{1}{N^\frac{1}{2}}].
\end{align}
From which one can get 
\begin{align}
	\sup_{t\in[0,N^{-\frac{1}{2}}]}\|\gamma(t)\|_{B_{\infty,1}^0}\geq CN^{\frac{3}{5}-\frac{1}{2}} \geq \ln N
\end{align}
which contradicts the hypothesis \eqref{as}.

 In conclusion, we obtain for $N>10$ large enough 
 \begin{align}
 	&\|\gamma(t)\|_{L_T^\infty( B_{\infty,1}^0)} \geq \ln N,~~~~T=N^{-\frac{1}{2}},
 	\\
 	&\|\gamma_0\|_{B_{\infty,1}^0} \leq N^{-\frac{1}{10}}
 \end{align}
that is we get the norm inflation and hence the ill-posedness of modified Camassa-Holm equation. Then we proved Theorem\ref{the1}

\end{proof}
\smallskip
\noindent\textbf{Acknowledgments} This work was
partially supported by the National Natural Science Foundation of China (No.11671407 and No.11701586), the Macao Science and Technology Development Fund (No. 098/2013/A3), and Guangdong Province of China Special Support Program (No. 8-2015),
and the key project of the Natural Science Foundation of Guangdong province (No. 2016A030311004).


\phantomsection
\addcontentsline{toc}{section}{\refname}
\bibliographystyle{abbrv} 
\bibliography{Feneref}

\end{document}